\documentclass[a4paper,10pt]{amsart}
\usepackage{amssymb,amsthm,amsmath,amscd,graphicx,enumerate,color,enumerate,mathtools}

\newtheorem{theorem}{Theorem}
\newtheorem*{theorem*}{Theorem}

\newtheorem{lemma}{Lemma}

\newtheorem{corollary}{Corollary}
\theoremstyle{definition}

\newcommand{\dd}{\ensuremath{\mathrm{d}}}
\newcommand{\ii}{\ensuremath{\mathrm{i}}}

\renewcommand{\Im}{\ensuremath{\mathsf{Im}\,}}
\renewcommand{\Re}{\ensuremath{\mathsf{Re}\,}}

\DeclarePairedDelimiter\floor{\lfloor}{\rfloor}

\title[Power mean of the Hurwitz zeta function]{Power Mean of the Hurwitz zeta function}

\begin{document}
\author[A.C.L. Ashton]{A.C.L. Ashton}
       \address{DAMTP, University of Cambridge.}
       \email{a.c.l.ashton@damtp.cam.ac.uk}

\date\today

\thanks{The author supported by Homerton College, University of Cambridge.}

\begin{abstract} 
In this note we derive asymptotic formulas for power mean of the Hurwitz zeta function over large intervals. 
\end{abstract}

\maketitle

\section{Introduction}
Throughout we use $s=\sigma+\ii t$ with $t>1$ and set $2\pi xy=t$ with $x\geq 1$ an integer. For $\sigma>1$ and $\alpha>0$ the Hurwitz zeta function is defined by
\[ \zeta(s,\alpha) = \sum_{m\geq 0} (m+\alpha)^{-s} \]
and by analytic continuation for $s\neq 1$. The modified Hurwitz zeta function will be defined as
\[ \zeta_x(s,\alpha) = \sum_{n\geq x} (n+\alpha)^{-s} = \zeta(s,\alpha) - \sum_{n<x} (n+\alpha)^{-s}, \]
with $\zeta_0(s,\alpha)=\zeta(s,\alpha)$. Here and throughout a sum over $n<x$ is to be interpreted as a sum over positive integers less than $x$. Much is known about the power means of $\zeta_1(s,\alpha)$ over the interval $(0,1)$. We investigate power means of $\zeta_1(u,\alpha)$ over a large $\alpha$-interval, i.e. the integrals
\[ \int_0^{t/2\pi}  |\zeta_1(\sigma+\ii t,\alpha)|^2\, \dd \alpha. \]
The asymptotic behaviour of $\zeta_1(s,\alpha)$ is very simple for $\alpha>t/2\pi$, so computation of power means over longer intervals would be very straightforward. Note that if $t/2\pi$ is an integer then this can be written
\begin{align*}
\int_0^{t/2\pi}  |\zeta_1(s,\alpha)|^2\, \dd \alpha &= \sum_{ x \leq t/2\pi} \int_{x-1}^x  |\zeta_1(s,\alpha)|^2\, \dd \alpha \\
&= \sum_{ x \leq t/2\pi} \int_{0}^1 |\zeta_x(s,\alpha)|^2\, \dd \alpha 
\end{align*}
and when $t/2\pi$ is not an integer the same expression holds upto an error term which is asymptotically smaller than the expression on the right hand side. This naturally leads us to the study of the integrals and sums of the form
\[ I_x(s) = \int_0^1 |\zeta_x(s,\alpha)|^2\, \dd \alpha, \quad \sum_{x\leq t/2\pi} I_x(s). \]

The integrals have already been the subject of much attention. For instance in \cite{rane1983} using a careful application of the Euler-Maclaurin theorem it was shown that
\begin{equation} I_x(s) = |K(s)|^2 \sum_{n\leq y} n^{2\sigma-2} + \mathcal{O}(x^{-2\sigma}) + \mathcal{O} \left( \frac{x^{2-2\sigma}}{t}\right)\label{Rane_estimate}\end{equation}
where $K(s)=(-2\pi \ii)^{s-1}\Gamma(1-s)$. This estimate is not useful when\footnote{We say $f\asymp g$ as $t\rightarrow \infty$ if $Af \leq g \leq Bf$ for some $B>A>0$ for $t$ sufficiently large.} $x\asymp t$ in the important case $\sigma=\tfrac{1}{2}$. Indeed, applying Stirling's estimate one finds
\[ |K\left( \tfrac{1}{2}+\ii t\right)|^2 \sum_{n\leq y} \frac{1}{n} = \log y + \gamma + o(1), \quad t\rightarrow \infty. \]
So when $x\asymp t$ the leading order term is $\log y=\mathcal{O}(1)$, which is the same size as the error term $\mathcal{O}(x/t)$. We require a sharper estimates to deal with the sums of interest, in which the range $x\asymp t$ is present. It transpires that sharper estimates required depend heavily on $\|y\|=\mathrm{dist}(y,\mathbf{Z})$.

The asymptotic expansion of $I_1(s)$ is well understood. A comprehensive derivation of results in this direction, based on Atkinson's disection \cite{atkinson1949}, can be found in \cite{katsurada1996explicit}. See also \cite{andersson1992,balasubramanian1979,rane1983} and references therein. In \cite{katsurada1996explicit} it is shown that for $\sigma\in(0,1)$ and $\sigma\neq \tfrac{1}{2}$
\begin{align*} 
I_1(s) &= \frac{1}{2\sigma-1}  + 2\Gamma(2\sigma -1) \zeta(2\sigma-1)\Re \left[ \frac{\Gamma(1-\sigma+\ii t)}{\Gamma(\sigma+\ii t)} \right]  \\
&\quad - 2\Re \sum_{n\geq 0}\frac{(\sigma+\ii t)_n }{(1-\sigma+\ii t)_{n+1}}\zeta_1(\sigma+\ii t+n,1). 
\end{align*}
Each term in the sum as asymptotically smaller than the preceding one, so this serves as a complete asymptotic expansion of $I_1(s)$ when $\sigma\neq \tfrac{1}{2}$. By setting $\sigma=\tfrac{1}{2}+\epsilon$ in this formula and taking the limit $\epsilon\rightarrow 0$ the authors obtain the remarkable asymptotic result
\[ I_1(\tfrac{1}{2}+\ii t) = \log \left( \frac{t}{2\pi}\right) +\gamma - 2 \Re \frac{\zeta(\tfrac{1}{2}+\ii t)}{\tfrac{1}{2}+\ii t} + \mathcal{O}\left( \frac{1}{t}\right). \]

Our aim is to provide a sharper estimate than in \eqref{Rane_estimate} valid for $\sigma \in (0,1)$. We do this in two ways, in section one we show that the order $\mathcal{O}(x^{2-2\sigma}/t)$ can be replaced with a sharper one that depends on an auxillary parameter $\eta\in (0,1/2)$. This error term vanishes if $x$ (recall $2\pi xy=t$) does not belong to the set
\[ A(t,\eta) = \left\{ 1\leq x \leq \frac{t}{2\pi}: \, \| y\| <\eta \right\} \]
where $\|y\| =\mathrm{dist}(y,\mathbf{Z})$ and $\eta \in (0,1/2)$. Using the results from \cite{saffari1977} we can estimate the size of $A(t,\eta)$ and by choosing $\eta=o(1)$ we show that the contribution from the sum over $0<x\leq t/2\pi$ for which $x\in A(t,\eta)$ is negligible. This gives the estimate
\[ \int_{0}^{t/2\pi} |\zeta_1(\sigma+\ii t,\alpha)|^2\, \dd \alpha = \left( \frac{t}{2\pi} \right)^{2-2\sigma} \zeta(3-2\sigma) + \mathcal{O} \left( t^{1-\sigma}\right) + \mathcal{O}\left( t^{7/4-2\sigma}\right), \]
valid for $\sigma\in (0,1)$ -- see Theorem \ref{thm3}.

\section{Statement of results I}
These first results are largely generalisations of those found in \cite{katsurada1996explicit}, for $x>1$. Throughout this section we refer to the set
\[ E= \{(u,v)\in \mathbf{C}^2: u+v\in \{2,1,0,-1,-2,\ldots\} \lor u\in \mathbf{Z} \lor v \in \mathbf{Z} \}, \]
and use the Pochhammer symbol $(s)_n = \Gamma(s+n)/\Gamma(s)$ for any integer $n$. We will use the notation
\[ J_x(u,v)=\int_0^1 \zeta_x(u,\alpha)\zeta_x(v,\alpha)\, \dd \alpha, \]
so that $I_s(s)=J_x(s,\bar{s})$. We write $A\lesssim B$ if $A\leq CB$ for some $C>0$. If the implied constant depends on a parameter $\beta$, we write $A\lesssim_\beta B$. The following result is a direct generalisation of Theorem 1 in \cite{katsurada1996explicit}, although our derivation is rather different. 
 
\begin{theorem}\label{lem1}
Let $N\geq 1$ be an integer, $u,v \in \mathbf{C}$ with $-N+1<\Re u < N+1$, $-N+1<\Re v < N+1$ and $(u,v)\notin E$. Then
\begin{align*} 
J_x(u,v) &= \frac{x^{1-u-v}}{u+v-1}  + \left[ \frac{\Gamma(1-u)}{\Gamma(v)} + \frac{\Gamma(1-v)}{\Gamma(u)}\right] \Gamma(u+v-1) \zeta(u+v-1) \\
& \quad - S_N(u,v;x) - S_N(v,u;x) - T_N(u,v;x) - T_N(v,u;x) 
\end{align*}
where 
\[ S_N(u,v;x) = \sum_{n=0}^{N-1} \frac{(u)_n x^{n+1-v}}{(1-v)_{n+1}} \zeta_x(u+n,1), \]
\[ T_N(u,v;x) = \frac{(u)_N x^{N+1-v}}{(1-v)_N} \sum_{l=1}^\infty l^{1-u-v} \int_l^\infty \beta^{u+v-2} (x+\beta)^{-u-N}\, \dd \beta. \]
In addition, for any $M\geq 0$
\begin{align*}
&\frac{T_N(u,v;x)}{x^{N+1-v}} = \sum_{m=1}^M (-1)^{m-1} \frac{(2-u-v)_{m-1}(u)_{N-m}}{(1-v)_N} \sum_{l=1}^\infty l^{-m} (x+l)^{-u-N+m} \\
&+ (-1)^M \frac{(2-u-v)_M (u)_{N-M}}{(1-v)_N} \sum_{l=1}^\infty l^{1-u-v} \int_l^\infty \beta^{u+v-M-2}(x+\beta)^{-u-N+M}\, \dd \beta.
\end{align*}
\end{theorem}

Following \cite{katsurada1996explicit} we note that of $N_0$ is chosen such that $-N_0+1<\Re u<N_0+1$ and $-N_0+1<\Re v<N_0+1$ for given $(u,v)\notin E$, then for $N\geq N_0$ we have
\begin{align*}
&\left|\sum_{l=1}^\infty l^{1-u-v} \int_l^\infty \beta^{u+v-2} (x+\beta)^{-u-N}\, \dd \beta\right| \\
&\quad \leq (x+1)^{-N+N_0} \sum_{l=1}^\infty l^{1-\Re (u+v)} \int_l^\infty \beta^{\Re (u+v)-2} (k+\beta)^{-\Re u - N_0}\, \dd \beta \\
&\quad \lesssim (x+1)^{-N},
\end{align*}
where the implied constant is independent of $N$. Using Stirling's approximation we then find that as $N\rightarrow \infty$
\[ T_N(u,v;x) = \mathcal{O} \left( N^{\Re (u+v)-1} \epsilon_x^N\right), \quad \epsilon_x = \frac{x}{x+1}<1. \]
Consequently, $\lim_{N\rightarrow \infty} T_N=0$ and we arrive at the first of our corollaries.
\begin{corollary}\label{cor1}
For $(u,v)\notin E$
\begin{align*} 
J_x(u,v) &= \frac{x^{1-u-v}}{u+v-1}  + \left[ \frac{\Gamma(1-u)}{\Gamma(v)} + \frac{\Gamma(1-v)}{\Gamma(u)}\right] \Gamma(u+v-1) \zeta(u+v-1) \\
&\quad - \sum_{n=0}^{\infty} \frac{(u)_n x^{n+1-v}}{(1-v)_{n+1}}  \zeta_x(u+n,1) - \sum_{n=0}^{\infty} \frac{(v)_n x^{n+1-u}}{(1-u)_{n+1}} \zeta_x(v+n,1). 
\end{align*}
\end{corollary}
When $u,v$ are complex conjugates, with $u=\sigma+\ii t$, we can get a simple estimate on $  T_N(u,v;x)$ using
\begin{align*}
&\left| \Re \sum_{l=1}^\infty l^{1-2\sigma} \int_l^\infty \beta^{2\sigma-M-2} (x+\beta)^{-\sigma-\ii t-N+M}\, \dd \beta \right| \\
&\qquad = \left| \frac{1}{t} \sum_{l=1}^{\infty} l^{1-2\sigma}  \int_l^\infty \frac{\beta^{2\sigma-M-2}}{ (x+\beta)^{\sigma+N-M-1}} \frac{\dd}{\dd \beta} \sin (t\log(x+\beta))\, \dd \beta \right| \\
&\qquad \lesssim \frac{1}{t} \sum_{l=1}^\infty \frac{1}{l^{M+1}(x+l)^{\sigma+N-M-1}} \\
&\qquad \lesssim \frac{x^{1-N-\sigma+M}}{t}
\end{align*}
valid for $M>0$. In the second line we applied the second mean value theorem for integrals and in the third we compared the sum to an appropriate integral. Similarly
\[ \left| \Im \sum_{l=1}^\infty l^{1-2\sigma} \int_l^\infty \beta^{2\sigma-2} (x+\beta)^{-\sigma - \ii t-N}\, \dd \beta \right| \lesssim \frac{x^{1-N-\sigma+M}}{t} \]
so on noting that
\[ \left| \frac{(2-2\sigma)_M(\sigma+\ii t)_{N-M}}{(1-\sigma+\ii t)_N} \right| \lesssim_{N,\sigma} t^{-M} \]
we arrive at the asymptotic estimate for $M>1$
\begin{align*}
T_N(u,v;x) &= \sum_{m=1}^M (-1)^{m-1} \frac{(2-u-v)_{m-1}(u)_{N-m}}{(1-v)_N} x^{N+1-v}\sum_{l=1}^\infty l^{-m} (x+l)^{-u-N+m} \\
&\qquad + \mathcal{O}_{N,M,\sigma}\left( \frac{x^{2-2\sigma +M}}{t^{M+1}}\right).
\end{align*}
This gives asymptotic expansion for $I_x(s)$ as $t\rightarrow \infty$.
\begin{corollary}\label{cor2}
For $N\geq 1$, $-N+1<\sigma<N+1$ and $2\sigma-1 \notin \{1,0,-1,-2,\ldots\}$
\begin{align*}
 I_x(s)&= \frac{x^{1-2\sigma}}{2\sigma-1} + 2 \Gamma(2\sigma-1) \zeta(2\sigma-1) \Re \left[ \frac{\Gamma(1-\sigma+\ii t)}{\Gamma(\sigma+\ii t)} \right] \\
&\quad - 2\Re \sum_{n=0}^{N-1} \frac{(\sigma+\ii t)_n x^{n+1-\sigma+\ii t}}{(1-\sigma+\ii t)_{n+1}} \zeta_x(\sigma+\ii t+n,1) +\mathcal{O}_{N,\sigma}\left( \frac{x^{2-2\sigma}}{t}\right).
\end{align*}
\end{corollary}
We see that the error term here is comparable to the leading order term when $x\asymp t$, in exactly the same way that the error term in \eqref{Rane_estimate} is. So these generalisations of the results of \cite{katsurada1996explicit} are not sufficient to get leading order estimates for the power mean of $\zeta_x(s,\alpha)$ over large intervals. However, this estimate is significantly sharper than \eqref{Rane_estimate} in the regime $x\asymp 1$.

To obtain the result in the important case on the critical line $\sigma=\frac{1}{2}$ we set $2\sigma=1+\epsilon$ and take a limit $\epsilon\rightarrow 0$, as in \cite{katsurada1996explicit}. In doing this we use the facts
\[ \frac{x^{1-2\sigma}}{2\sigma-1} = \frac{1}{\epsilon} - \log x + o(1), \quad \Gamma(2\sigma-1) = \frac{1}{\epsilon} - \gamma + o(1) \]
as $\epsilon\rightarrow 0$. Using these results in the asymptotic estimate for $I_x(s)$ with $N=1$ we obtain the third of our corollaries.

\begin{corollary}\label{cor3}
For $x\geq 1$ an integer and $2\pi xy=t$
\[ I_x(\tfrac{1}{2}+\ii t)= \log y+ \gamma - 2 \Re \frac{x^{\frac{1}{2}+\ii t}\zeta_x(\tfrac{1}{2}+\ii t,1)}{\tfrac{1}{2}+\ii t} +  \mathcal{O}\left( \frac{x}{t} \right). \]
\end{corollary}

\section{Statement of results II}
The results in this section are found by calculating $I_x(s)$ in a different way that allows for a more precise remainder term. This form of remainder allows for the calculation of the power mean of $\zeta_1(s,\alpha)$ over a large interval. We remind the reader of the set
\[ A(t,\eta) = \left\{ 1\leq x \leq \frac{t}{2\pi}: \, \| y\| <\eta \right\}. \]
\begin{theorem}\label{thm2}
Let $s=\sigma+\ii t$, $t>1$, $2\pi xy=t$ with $x\geq 1$ an integer, $\sigma \in (0,1)$ and $\eta\in (0,1/2)$. Then for $x\notin A(t,\eta)$
\[ I_x(s) = |K(s)|^2 \sum_{m\leq y-\eta} m^{2\sigma-2} + \mathcal{O}\left( \frac{x^{-2\sigma}}{\eta^2}\right) + \mathcal{O}\left( \frac{t^{-1/2} x^{1-2\sigma}\log(y+2)}{\eta}\right)  \]
and for $x\in A(t,\eta)$ 
\begin{multline*}
I_x(s) = |K(s)|^2 \sum_{m\leq y-\eta} m^{2\sigma-2} + \mathcal{O}\left( \frac{x^{-2\sigma}}{\eta^2}\right) + \mathcal{O}\left( \frac{t^{-1/2} x^{1-2\sigma}\log(y+2)}{\eta}\right) \\
+ \left( \frac{t}{2\pi}\right)^{1-2\sigma} [y]^{2\sigma-2}  \left|\mathcal{E}\left( t, \frac{[y]-y}{y}\right)\right|^2 + \mathcal{O}\left( \frac{x^{2-2\sigma}}{t^{3/2}}\right).
\end{multline*}
Here $[y]$ denotes the nearest integer to $y$, $\mathcal{E}(t,0)=1/2$ and
\[ \mathcal{E}(t,a)=  H(-a) + e^{\ii t a-\ii \pi/4}\mathrm{sgn}(a)(1+a)^{-it}\int_0^{\infty} e^{\ii \pi \tau^2} e^{2\pi \ii\tau \sqrt{\frac{t}{2\pi}} |a|}\, \dd \tau , \]
with $H$ denoting the Heaviside function.
\end{theorem}
Note that there is no ambiguity in $[y]$ if $x\in A(t,\eta)$, since $\|y\|<\eta<1/2$.

\begin{theorem}\label{thm3}
For each $\sigma\in(0,1)$ and $t>1$
\[ \int_{0}^{t/2\pi} |\zeta_1(\sigma+\ii t,\alpha)|^2\, \dd \alpha = \left( \frac{t}{2\pi} \right)^{2-2\sigma} \zeta(3-2\sigma) + \mathcal{O} \left( t^{1-\sigma}\right) + \mathcal{O}\left( t^{7/4-2\sigma}\right). \]
\end{theorem}

It is important to reconcile the result in Theorem \ref{thm2} with that of Corollary \ref{cor2} in the previous section. We note from the functional equation
\[ 2 \Gamma(2\sigma-1)\zeta(2\sigma-1) =  \frac{\zeta(2-2\sigma)}{(2\pi)^{1-2\sigma}  \sin(\pi \sigma)} \]
so Stirling's approximation gives
\[ 2 \Gamma(2\sigma-1) \zeta(2\sigma-1) \Re \left[ \frac{\Gamma(1-\sigma+\ii t)}{\Gamma(\sigma+\ii t)} \right] =  \left( \frac{t}{2\pi}\right)^{1-2\sigma} \zeta(2-2\sigma) + \mathcal{O}\left( t^{-1-2\sigma}\right). \]
The expression appearing in Corollary \ref{cor2} can now be written as
\begin{align*}
 I_x(s)&= \left(\frac{t}{2\pi}\right)^{1-2\sigma} \left[ \frac{y^{1-2\sigma}}{2\sigma-1} + \zeta(2-2\sigma)\right] \\
&\quad - 2\Re \sum_{n=0}^{N-1} \frac{(\sigma+\ii t)_n x^{n+1-\sigma+\ii t}}{(1-\sigma+\ii t)_{n+1}} \zeta_x(\sigma+\ii t+n,1) +\mathcal{O}_{N,\sigma}\left( \frac{x^{2-2\sigma}}{t}\right),
\end{align*}
where we have absorbed the $\mathcal{O}(t^{-1-2\sigma})$ into the final error term. Now for $\sigma<1/2$ we may write
\[ \zeta(2-2\sigma)- \sum_{m\leq y-\eta} m^{2\sigma-2} = \sum_{m>y-\eta} m^{2\sigma-2}. \]
and by Euler-Maclaurin the following holds
\[ \sum_{m>y-\eta}  m^{2\sigma-2} =  - \frac{(y-\eta)^{2\sigma-1}}{2\sigma-1}+ \left( \{y-\eta\}-\tfrac{1}{2}\right) (y-\eta)^{2\sigma-2} + \int_{y-\eta}^\infty \frac{ \{\alpha\}-\tfrac{1}{2}}{ \alpha^{3-2\sigma}}\, \dd \alpha. \]
The right hand provides an analytic continuation of the left hand side for $\sigma<1$. Estimating the integral using the second mean value theorem for integrals and expanding the first term binomially, recalling $\eta\in(0,1/2)$ and $y\geq 1$, we find
\[ \zeta(2-2\sigma)- \sum_{m\leq y-\eta} m^{2\sigma-2} = - \frac{y^{2\sigma-1}}{2\sigma-1} + \left( \{y-\eta\}+\eta - \tfrac{1}{2}\right) y^{2\sigma-2} + \mathcal{O}(y^{2\sigma-3}). \]
Using $|K(s)|^2=t^{1-2\sigma} + \mathcal{O}(t^{-1-2\sigma})$ we conclude
\begin{align*}
 I_x(s)&=  |K(s)|^2 \sum_{m\leq y-\eta} m^{2\sigma-2} \\
&\quad - 2\Re \sum_{n=0}^{N-1} \frac{(\sigma+\ii t)_n x^{n+1-\sigma+\ii t}}{(1-\sigma+\ii t)_{n+1}} \zeta_x(\sigma+\ii t+n,1) +\mathcal{O}_{N,\sigma}\left( \frac{x^{2-2\sigma}}{t}\right),
\end{align*}
where we absorbed the $\mathcal{O}(t^{1-2\sigma} y^{2\sigma-2})$ term into the $\mathcal{O}(x^{2-2\sigma}/t)$ term. So this result is consistent with Corollary \ref{cor2}, but clearly the estimates in Theorem \ref{thm2} are sharper when $x\asymp t$.

\section{Proof to Theorem 1}
Fix an integer $N\geq 1$. Let us first assume that $1<\Re u <N+1$ and $1<\Re v <N+1$. The more general result will follow from a simple analytic continuation argument. In this case we have the identity
\[ \frac{1}{(m+\alpha)^{u}} = \frac{1}{\Gamma(u)} \int_0^\infty p^{u-1} e^{-p(m+\alpha)}\, \dd p, \]
where $\Gamma(u)$ denotes the Gamma function. Using this we find
\[ \zeta_x(u,\alpha) = \frac{1}{\Gamma(u)} \int_0^\infty \frac{p^{u-1} e^{-\alpha p} e^{-(x-1)p}}{e^p-1}\, \dd p \]
where $x\geq 1$. Noting that the function
\[ (p,q,\alpha;u,v)\mapsto \frac{p^{\Re(u)-1} e^{-\alpha p} e^{-(x-1)p}}{e^p-1} \frac{q^{\Re (v)-1} e^{-\alpha q} e^{-(x-1)q}}{e^q-1} \]
is absolute integrable $(0,\infty)\times (0,\infty)\times (0,1)$ for $\Re u, \Re v>1$ we can apply Fubini's theorem to interchange the orders of integration so that 
\begin{align*}
J_x(u,v) &= \frac{1}{\Gamma(u)\Gamma(v)} \int_0^\infty \int_0^\infty  \left[ \frac{p^{u-1} q^{v-1} e^{-(x-1)(p+q)}}{(e^p-1)(e^q-1)}\right] \left[ \frac{1-e^{-(p+q)}}{p+q}\right]\, \dd p\, \dd q \\
& = \frac{1}{\Gamma(u)\Gamma(v)} \lim_{\epsilon\downarrow 0} \iint_{D_\epsilon}  \left[ \frac{p^{u-1} q^{v-1} e^{-(x-1)(p+q)}}{(e^p-1)(e^q-1)}\right] \left[ \frac{1-e^{-(p+q)}}{p+q}\right]\, \dd p\, \dd q
\end{align*}
where $D_\epsilon=\{(p,q):p>\epsilon,q>\epsilon\}$. Following \cite{ashtonfokas2017} we observe the identity
\[ \frac{1}{(e^p-1)(e^q-1)}\frac{1-e^{-(p+q)}}{p+q} \equiv \left[ \frac{1}{e^p-1} + \frac{1}{e^q-1}+1 \right] \frac{e^{-(p+q)}}{p+q}. \]
Using this in the previous expression we arrive at
\begin{align}
J_x(u,v) &= r_x(u,v)+r_x(v,u) \label{Ik_eq1} \\
&\quad +\frac{1}{\Gamma(u)\Gamma(v)} \lim_{\epsilon\downarrow 0}\iint_{D_\epsilon} p^{u-1} q^{v-1} e^{-(x-1)(p+q)} \left[ \frac{ e^{-(p+q)}}{p+q}\right] \dd p\, \dd q \nonumber
\end{align}
where we have defined
\[ r_x(u,v) = \frac{1}{\Gamma(u)\Gamma(v)} \lim_{\epsilon\downarrow 0}\iint_{D_\epsilon} \left[ \frac{p^{u-1} e^{-p(x-1)}}{e^p-1} \right] \left[ \frac{e^{-(p+q)}}{p+q}\right] q^{v-1} e^{-q(x-1)}\, \dd p\, \dd q. \]
For each $\epsilon>0$ and $(p,q)\in D_\epsilon$ we have the identity
\[ \frac{e^{-(p+q)}}{p+q} = \int_1^\infty e^{-\alpha(p+q)}\, \dd \alpha. \]
Again, owing to the absolute convergence of the relevant triple integrals, we can interchange orders of integration so that
\begin{align*} &\frac{1}{\Gamma(u)\Gamma(v)} \lim_{\epsilon\downarrow 0}\iint_{D_\epsilon} p^{u-1} q^{v-1} e^{-(x-1)(p+q)} \left[ \frac{ e^{-(p+q)}}{p+q}\right] \dd p\, \dd q \\
&\quad = \frac{1}{\Gamma(u)\Gamma(v)}  \lim_{\epsilon\downarrow 0}\int_1^\infty \iint_{D_\epsilon}  p^{u-1} q^{v-1} e^{-(\alpha+x-1)(p+q)}  \dd p\, \dd q\, \dd \alpha \\
&\quad = \int_1^\infty (\alpha+x-1)^{-u-v}\, \dd \alpha \\
&\quad = \frac{x^{1-u-v}}{u+v-1}
\end{align*}
and similarly
\begin{align*} 
r_x(u,v) &= \frac{1}{\Gamma(u)\Gamma(v)} \lim_{\epsilon\downarrow 0}\int_1^\infty \iint_{D_\epsilon}  \left[ \frac{p^{u-1} e^{-p(\alpha+ x-1)}}{e^p-1} \right]q^{v-1} e^{-q(\alpha+x-1)}\, \dd p\, \dd q\, \dd\alpha \\
& =\int_1^\infty (\alpha+x-1)^{-v} \zeta_x(u,\alpha)\, \dd \alpha \\
& =\int_x^\infty \alpha^{-v} \zeta_1(u,\alpha)\, \dd \alpha.
\end{align*}
Using both these results \eqref{Ik_eq1} becomes
\begin{equation}J_x(u,v) = \frac{x^{1-u-v}}{u+v-1} + \int_x^\infty \alpha^{-u}\zeta_1(v,\alpha)\,\dd\alpha + \int_x^\infty \alpha^{-v}\zeta_1(u,\alpha)\,\dd\alpha. \label{Ik_eq2} \end{equation}
Recall the estimate $\zeta_1(u,\alpha)=\mathcal{O}(\alpha^{1-u})$ as $\alpha\rightarrow \infty$, see \cite{ferreira2004asymptotic} and references therein. It is clear then that both integrands are $\mathcal{O}(\alpha^{1-u-v})$ for large $\alpha$, so the corresponding integrals converge absolutely when $\Re u,\Re v>1$. Integrating by parts $N$ times and repeatedly applying the formula $\partial_\alpha \zeta_1(u,\alpha)=-u\zeta_1(u+1,\alpha)$ we find
\[ \int_x^\infty \alpha^{-v}\zeta_1(u,\alpha)\,\dd \alpha  = -S_N(u,v;x) + \frac{(u)_N}{(1-v)_N} \int_x^\infty \alpha^{-v+N}\zeta_1(u+N,\alpha)\,\dd\alpha. \]
Again following \cite{ashtonfokas2017} we have that for $\Re b<1$ and $\Re (a+b)>2$
\[ \int_0^\infty \alpha^{-b} \zeta_1(a,\alpha)\, \dd \alpha = \frac{\Gamma(1-b)}{\Gamma(a)} \Gamma(a+b-1) \zeta(a+b-1). \]
If $1<\Re u<N+1$ and $1<\Re v<N+1$ then this result is valid with the choices $a=u+N$ and $b=v-N$, so we find
\begin{align*} & \int_x^\infty \alpha^{-v+N}\zeta_1(u+N,\alpha)\,\dd\alpha \\
&\quad = \frac{\Gamma(1-v+N)}{\Gamma(u+N)} \Gamma(u+v-1) \zeta(u+v-1) - \int_0^x \alpha^{-u+N} \zeta_1(u+N,\alpha)\,\dd\alpha.
\end{align*}
Using these results in \eqref{Ik_eq2} we find
\begin{align}
J_x(u,v) &= \left[ \frac{\Gamma(1-u)}{\Gamma(v)}+ \frac{\Gamma(1-v)}{\Gamma(u)}\right] \Gamma(u+v-1)\zeta(u+v-1) \label{Ik_eq3}\\
&\quad -S_N(u,v;x) - S_N(v,u;x)  \nonumber \\
&\quad - \frac{(u)_N}{(1-v)_N} \int_0^x \alpha^{-v+N} \zeta_1(u+N,\alpha)\,\dd\alpha \nonumber  \\
&\quad -  \frac{(v)_N}{(1-u)_N} \int_0^x \alpha^{-u+N} \zeta_1(v+N,\alpha)\,\dd\alpha \nonumber  
\end{align}
For the final terms we write
\[ \int_0^x \alpha^{-v+N} \zeta_1(u+N,\alpha)\,\dd\alpha = \sum_{l=1}^\infty \int_0^x \alpha^{-v+N} (\alpha+l)^{-u-N}\, \dd \alpha, \]
the interchange being justified by the uniform convergence of the sum. Making the substitution $\alpha = lx/\beta$ gives the first form of $T_N(u,v;x)$ and the second form comes from a simple integration by parts $M$ times. The result in the extended ranges of $\Re u, \Re v$ follows by analytic continuation.

\section{Proof of Theorem 2}
In this section we derive estimates for $I_x(s)$ using a similar approach to \cite{rane1997}. Using Euler-Maclaurin it can be shown that for $\sigma>0$
\[ \zeta_1(s,\alpha) = \frac{\alpha^{1-s}}{s-1} - \frac{\alpha^{-s}}{2} + \sum_{m\neq 0} \left( \int_{\alpha}^\infty \beta^{-s} e^{2\pi\ii m \beta}\, \dd \beta\right) e^{-2\pi\ii m \alpha}, \]
where the notation $\sum_{m\neq 0}$ implicitly means $\lim_{N\rightarrow \infty} \sum_{0<|m|<N}$. We note from \cite{ashtonfokas2017} that for $\alpha>t/2\pi+\delta$ for $\delta\in(0,1)$ we have
\[ \left| \sum_{0<|m|<N} \left( \int_{\alpha}^\infty \beta^{-s} e^{2\pi\ii m \beta}\, \dd \beta\right) e^{-2\pi\ii m \alpha} \right| \lesssim_\delta t \alpha^{-\sigma-1}, \]
uniformly in $N$, so that in particular, the function
\[ \alpha\mapsto \zeta_1(s,\alpha) - \frac{\alpha^{1-s}}{s-1} + \frac{\alpha^{-s}}{2} \]
is absolutely integrable for $\sigma>0$. We rewrite \eqref{Ik_eq2} as
\begin{align*} J_x(u,v) &= \frac{x^{1-u-v}}{u+v-1} + \int_x^\infty \alpha^{-u}\zeta_1(v,\alpha)\,\dd\alpha + \int_x^\infty \alpha^{-v}\zeta_1(u,\alpha)\,\dd\alpha. 
 \\
& = \frac{1}{2}\frac{x^{1-u-v}}{u+v-1} + \int_x^\infty \alpha^{-u} \left(\frac{\alpha^{1-v}}{v-1} - \frac{\alpha^{-v}}{2} \right)\, \dd \alpha \\
&\quad + \sum_{m\neq 0} \int_x^\infty  \left( \int_{\alpha}^\infty \beta^{-v} e^{2\pi\ii m \beta}\, \dd \beta\right) \alpha^{-u} e^{-2\pi\ii m \alpha}\, \dd\alpha + (u\leftrightarrow v)
\end{align*}
Performing the integrals and re-ordering one of the sums we find
\begin{align*} J_x(u,v) &= \frac{x^{2-u-v}}{(u-1)(v-1)}  + \sum_{m\neq 0} \int_x^\infty  \left( \int_{\alpha}^\infty \beta^{-u} e^{2\pi\ii m \beta}\, \dd \beta\right) \alpha^{-v} e^{-2\pi\ii m \alpha}\, \dd \alpha \\
&\quad + \sum_{m\neq 0} \int_x^\infty  \left( \int_{\alpha}^\infty \beta^{-v} e^{-2\pi\ii m \beta}\, \dd x\right) \alpha^{-u} e^{2\pi\ii m \alpha}\, \dd \alpha.
\end{align*}
The $\alpha$-integrands in the latter terms can be written
\[ - \frac{\dd}{\dd \alpha}\left[ \left( \int_{\alpha}^\infty \beta^{-u} e^{2\pi\ii m \beta}\, \dd \beta\right)\left( \int_{\alpha}^\infty \beta^{-v} e^{-2\pi\ii m \beta}\, \dd \beta\right)\right] \]
so we find
\begin{equation}
J_x(u,v) = \frac{x^{2-u-v}}{(u-1)(v-1)} + \sum_{m\neq 0} \left( \int_{x}^\infty \beta^{-u} e^{2\pi\ii m \beta}\, \dd \beta\right)\left( \int_{x}^\infty \beta^{-v} e^{-2\pi\ii m \beta}\, \dd \beta\right)
\end{equation}
and for the particular case $u=s=\bar{v}$ we find
\begin{equation} I_x(s) = \frac{x^{2-2\sigma}}{t^2 + (\sigma-1)^2} + \sum_{m\neq 0} \left| \int_x^\infty \beta^{-s} e^{2\pi \ii m \beta}\, \dd \beta\right|^2. \label{Fourier} \end{equation}
Note a similar result in the case $x=1$ was obtained by an entirely different method in \cite{rane1997}. The integral in the sum is easy to estimate if either $m<0$ or $m>y + \eta$ for $\eta \in (0,1/2)$. Indeed, we see that
\begin{align*} \int_x^\infty \beta^{-s} e^{2\pi \ii m \beta}\, \dd \beta &= \frac{1}{2\pi \ii }\int_x^\infty \frac{ \beta^{-\sigma}}{m - \frac{t}{2\pi \beta}} \frac{\dd}{\dd \beta} \left( e^{\ii t(2\pi m \beta/t - \log \beta)}\right) \dd \beta \\
&\equiv \frac{1}{2\pi \ii }\int_x^\infty F(\beta) \frac{\dd}{\dd \beta} \left( e^{\ii t(2\pi m \beta/t - \log \beta)}\right) \dd \beta.
\end{align*}
Noting that $F$ $(-F)$ is positive and monotone decreasing for $m>y+\eta$ $(m<0)$, employing the second mean value theorem for integrals on the real and imaginary parts of this integral we deduce that
\[ \sum_{m<0} \left| \int_x^\infty \beta^{-s} e^{2\pi \ii m \beta}\, \dd \beta \right|^2 =\mathcal{O}\left( \sum_{m<0} \frac{x^{-2\sigma}}{(m-y)^2}\right) = \mathcal{O}(x^{-2\sigma}) \]
and similarly
\[ \sum_{m>y + \eta} \left| \int_x^\infty \beta^{-s} e^{2\pi \ii m \beta}\, \dd \beta \right|^2 = \mathcal{O}\left( \sum_{m>y+\eta}\frac{x^{-2\sigma}}{(m-y)^2}\right)=\mathcal{O}\left( \frac{x^{-2\sigma}}{\eta^2}\right). \]
For $1\leq m < y -\eta$ we write
\begin{align*}
 \int_x^\infty \beta^{-s} e^{2\pi\ii m \beta}\, \dd \beta &= \int_0^\infty \beta^{-s} e^{2\pi\ii m \beta}\, \dd \beta - \int_0^x \beta^{-s} e^{2\pi\ii m\beta}\, \dd \beta \\
 &= K(s) m^{s-1} - \frac{1}{2\pi\ii}\int_0^x \frac{\beta^{-\sigma}}{m-\frac{t}{2\pi\beta}} \frac{\dd}{\dd \beta}\left( e^{\ii t(2\pi m \beta/t - \log \beta)}\right)\dd \beta 
\end{align*}
where $K(s)=(2\pi/\ii)^{s-1}\Gamma(1-s)$. Stirling's approximation gives
\[ |K(s)|^2 = \left( \frac{t}{2\pi}\right)^{1-2\sigma} + \mathcal{O}\left( t^{-2\sigma}\right) \quad \textrm{as $t\rightarrow$} \infty \]
so again employing the second mean value for integrals we find
\begin{align*} & \sum_{m<y -\eta} \left| \int_x^\infty \beta^{-s} e^{2\pi\ii m \beta}\, \dd \beta\right|^2 \\ & \qquad = |K(s)|^2 \sum_{m<y-\eta} m^{2\sigma-2} + \sum_{0<m<y-\eta}\mathcal{O}\left( |K(s)|\frac{ m^{\sigma-1} x^{-\sigma}}{y-m} + \frac{x^{-2\sigma}}{(y-m)^2}\right) \\
& \qquad  = |K(s)|^2  \sum_{m<y-\eta} m^{2\sigma-2} + \mathcal{O}\left( \frac{t^{1/2-\sigma} y^{\sigma-1} x^{-\sigma}\log (y+2)}{\eta}\right) + \mathcal{O}\left( \frac{x^{-2\sigma}}{\eta^2}\right) \\
& \qquad  = \left( \frac{t}{2\pi}\right)^{1-2\sigma}  \sum_{m<y-\eta} m^{2\sigma-2} + \mathcal{O}\left( \frac{t^{-1/2} x^{1-2\sigma}\log (y+2)}{\eta}\right) + \mathcal{O}\left( \frac{x^{-2\sigma}}{\eta^2}\right). 
\end{align*}
In the final line we absorbed an error term of $\mathcal{O}\left( t^{-2\sigma} y^{2\sigma-1}\right)=\mathcal{O}(t^{-1} x^{-2\sigma})$ into the third error term appearing on the right hand side. In the second line we used the simple estimate
\[ \sum_{m<y-\eta} \frac{m^{\sigma-1}}{y-m} = \mathcal{O}\left(\frac{y^{\sigma-1}}{\eta}\right) + \sum_{m<y-\eta-1} \frac{m^{\sigma-1}}{y-m} \]
and the latter sum can be estimated by first splitting into ranges in which the summand is increasing/decreasing then estimating each part by comparing to an appropriate integral. We conclude
\[ \sum_{m<y-\eta} \frac{m^{\sigma-1}}{y-m} = \mathcal{O}\left( \frac{y^{\sigma-1}}{\eta}\right) + \mathcal{O}\left( y^{\sigma-1} \log (y+2) \right) = \mathcal{O}\left( \frac{y^{\sigma-1} \log (y+2)}{\eta}\right). \]

To estimate the sum on the (possibly empty) range $|m-y|\leq \eta$ more care is needed because of the stationary point at $\beta=t/2\pi m$. It will be convenient to make a change of variables that fixes the stationary point at the origin. To this end we set
\[ \beta'= \left(\frac{2\pi m}{t}\right) \beta -1 \]
which gives, after dropping the primes,
\begin{equation} \int_x^\infty \beta^{-s} e^{2\pi \ii m \beta}\, \dd \beta  = e^{\ii t}\left( \frac{t}{2\pi m}\right)^{1-s}  \int_{a}^\infty (\beta+1)^{-\sigma} e^{\ii t(\beta - \log (\beta+1))}\, \dd \beta, \label{int-estimator} \end{equation}
where we have defined
\[ a =  \frac{m- y}{y} . \]
Note that $|a|<1$ for $|m-y|\leq\eta$, and $|a|$ can come very close to zero, which is the location of the stationary point. For this reason we need an estimate which is uniform with respect to the distance of the stationary point from the end-point of integration. Uniform asymptotics of integrals of this kind are dealt with in \cite{skinner1997} (c.f. also \cite{bleistein1966}), and the relevant form of the appropriate lemma is as follows.

\begin{lemma}\label{asymplem}
Fix $b>0$. For $g,h \in C^\infty[0,b]$ with $g(0)=1$, $h(0)=h'(0)=0$, $h''(0)=1$ and $h'(x)>0$ for $0<\beta\leq b$ it holds that
\[ \int_0^a g(\beta) e^{\ii th(\beta)}\, \dd \beta = \frac{1}{2} e^{\ii \pi/4} \sqrt{\frac{2\pi }{t }} \left[ 1  - e^{\ii t h(a)}\Psi\left( \sqrt{\frac{t }{2\pi}} a\right) \right] + \mathcal{O}\left(\frac{1}{t}\right)  \]
as $t\rightarrow \infty$, uniformly in $0\leq a \leq b$. Here $\Psi$ denotes the Fresnel-type integral
\[ \Psi(\epsilon) = \frac{2 e^{-\ii \pi \epsilon^2}}{ e^{\ii \pi/4}}  \int_\epsilon^{\infty} e^{\ii \pi s^2}\, \dd s \equiv 2 e^{-\ii \pi/4} \int_0^\infty e^{\ii \pi s^2} e^{2\pi\ii s \epsilon}\, \dd s.\]
\end{lemma}

In the case $a>0$ the the integral appearing on the right hand side of \eqref{int-estimator} is
\[ \int_a^\infty g(\beta) e^{\ii th(\beta)}\, \dd \beta, \qquad g(\beta)=(1+\beta)^{-\sigma}, \quad h(\beta) = \beta-\log (1+\beta), \]
which we write as
\[\int_a^\infty g(\beta) e^{\ii th(\beta)}\, \dd \beta =  \int_0^\infty g(\beta)e^{\ii th(\beta)}\, \dd \beta - \int_0^a g(\beta) e^{\ii th(\beta)}\, \dd\beta. \]
The first integral is gives (half) a stationary phase contribution from the stationary point at $\beta=0$ with a uniform error of $\mathcal{O}(1/t)$. The second integral can be estimated using Lemma \ref{asymplem}. We find that for $a>0$, i.e. $m>y$,
\begin{align}& \int_x^\infty \beta^{-s} e^{2\pi\ii m \beta}\, \dd \beta \label{SPapos}\\
&\qquad  = e^{\ii t+\ii \pi/4} \left( \frac{t}{2\pi m}\right)^{1-s} \sqrt{\frac{2\pi}{t}}\left[ \frac{e^{\ii th(a)}}{2} \Psi\left( \sqrt{\frac{t}{2\pi}} a\right)\right] + \mathcal{O}\left( \frac{m^{\sigma-1}}{t^\sigma}\right)\nonumber .
\end{align}
If $a<0$ we write the integral of the right hand side of \eqref{int-estimator} as
\begin{align*}
\int_a^\infty g(\beta) e^{\ii th(\beta)}\, \dd \beta &=  \int_0^\infty g(\beta)e^{\ii th(\beta)}\, \dd \beta + \int_a^0 g(\beta) e^{\ii th(\beta)}\, \dd\beta\\
&= \int_0^\infty g(\beta)e^{\ii th(\beta)}\, \dd \beta + \int_0^{|a|} \bar{g}(x) e^{\ii t \bar{h}(t)}\, \dd \beta,
\end{align*}
where $\bar{g}(\beta) = (1-\beta)^{-\sigma}$ and $\bar{h}(\beta) = -\beta-\log(1-\beta)$. Again $\bar{g}$ and $\bar{h}$ satisfy the conditions of Lemma \ref{asymplem} so arguing as before we find for $a<0$, i.e. $m<y$
\begin{align} & \int_x^\infty \beta^{-s} e^{2\pi\ii m \beta}\, \dd \beta \label{SPaneg}\\
&\qquad = e^{\ii t+\ii \pi/4} \left( \frac{t}{2\pi m}\right)^{1-s}\sqrt{\frac{2\pi}{t}} \left[ 1- \frac{e^{\ii t h(a)}}{2} \Psi\left( \sqrt{\frac{t}{2\pi}} |a|\right)\right]+ \mathcal{O}\left( \frac{m^{\sigma-1}}{t^\sigma}\right). \nonumber
\end{align}
Note that these results are consistent with in the limit $a\downarrow 0$, where $m=y$. In conclusion, for $|m-y|<\eta$
\begin{equation} \int_x^\infty \beta^{-s} e^{2\pi\ii m \beta}\, \dd \beta = e^{\ii t+\ii \pi/4} \left( \frac{t}{2\pi m}\right)^{1-s}\sqrt{\frac{2\pi}{t}} \mathcal{E}(t,a)+ \mathcal{O}\left( \frac{m^{\sigma-1}}{t^\sigma}\right), \label{unif_int} \end{equation}
where $\mathcal{E}(t,0)=1/2$ and for $a\neq 0$
\[ \mathcal{E}(t,a)=  H(-a) + \mathrm{sgn}(a)\frac{e^{\ii t h(a)}}{2} \Psi\left( \sqrt{\frac{t}{2\pi}} |a|\right), \]
where $H$ denotes the Heaviside function. We see from the definitions that $\mathcal{E}=\mathcal{O}(1)$. Noting that $|m-y|=\|y\|$ if $|m-y|<\eta<1/2$, and $m=[y]$, where $[y]$ denotes the nearest integer to $y$, the right hand side of \eqref{unif_int} can be written
\[ \frac{e^{\ii t+\ii \pi/4}}{2} \left( \frac{t}{2\pi [y]}\right)^{1-s} \sqrt{\frac{2\pi}{t}} \mathcal{E}\left( t, \frac{[y]-y}{y}\right) + \mathcal{O}\left( \frac{x^{1-\sigma}}{t} \right). \]
In summary
\[ \sum_{|m-y|\leq \eta} \left| \int_x^\infty \beta^{-s} e^{2\pi\ii m \beta}\, \dd \beta\right|^2 = \frac{1}{4} \left( \frac{t}{2\pi}\right)^{1-2\sigma} [y]^{2\sigma-2}  \left|\mathcal{E}\left( t, \frac{[y]-y}{y}\right)\right|^2 + \mathcal{O}\left( \frac{x^{2-2\sigma}}{t^{3/2}}\right). \]
Note that this term contributes \emph{nothing} unless $\|y\|<\eta$, i.e. $x$ belongs to $A(t,\eta)$. Collecting the different cases $m<0$, $m>y+\eta$, $m<y-\eta$ and $|m-y|\leq \eta$ we arrive at the result in Theorem \ref{thm2}.

\section{Proof of Theorem 3}
We will be required to estimate the size of the set
\[ A(t,\eta) = \left\{ 1\leq x \leq \frac{t}{2\pi}: \, \| y\| <\eta \right\}. \]
The following lemma can be found in Corollary 1.1 of \cite{saffari1977}.
\begin{lemma}
Let $\{x\}$ denote the fractional part of $x$. Then for each $\delta\in (0,1)$
\[ \frac{\# \{ 1\leq x \leq n: \{n/x\}<\delta \}}{n} = \sum_{m\geq 1} \frac{\delta }{m(m+\delta)} + \mathcal{O}\left( n^{-2/3} \log n\right)\,\, \textrm{as $n\rightarrow \infty$}. \]
The implicit constant in the error term is uniform in $\delta$.
\end{lemma}
As the authors state in \cite{saffari1977}, it is likely that the error term can be improved using more sophisticated machinery. However, the error given here is sufficient for our purposes. Using the digamma function $\psi=\Gamma'/\Gamma$ this result can be restated as
\begin{equation} \# \{ 1\leq x \leq n: \{n/x\}<\delta \} = \left( \gamma + \psi(1+\delta)\right) n + \mathcal{O}\left( n^{1/3} \log n\right). \label{meas1} \end{equation}
By taking compliments we also have the result
\begin{equation} \# \{ 1\leq x \leq n: \{n/x\}>1-\delta \} = n - \left( \gamma + \psi(2-\delta) \right)n + \mathcal{O}\left( n^{1/3} \log n\right). \label{meas2} \end{equation}
It is clear from the definitions that since $\eta\in (0,1/2)$
\[ A(t,\eta) = \left\{ 1\leq x \leq \tfrac{t}{2\pi}: \{t/2\pi x\}<\eta \right\} \cup \left\{ 1\leq x \leq \tfrac{t}{2\pi}: \{t/2\pi x\}>1-\eta \right\},\] 
and each of these sets are disjoint. On using estimates \eqref{meas1} and \eqref{meas2} with $n=t/2\pi$ we deduce
\begin{align*} \# A(t,\eta) &= \frac{t}{2\pi}\left( 1+ \psi(1+\eta) - \psi(2-\eta)\right) + \mathcal{O}(t^{1/3}\log t)\\
 &= \mathcal{O}(\eta t) + \mathcal{O}(t^{1/3}\log t). 
\end{align*}
We perform the sum
\[ \sum_{ x \leq t/2\pi} I_s(s) = \sum_{A^c} I_x(s) + \sum_{A}I_x(s). \]
From Theorem \ref{thm2} and our previous result we see that the right hand side is
\begin{multline*}
\sum_{ x \leq t/2\pi} \left( |K(s)|^2 \sum_{m\leq y-\eta} m^{2\sigma-2} + \mathcal{O}\left( \frac{x^{-2\sigma}}{\eta^2}\right) + \mathcal{O}\left( \frac{t^{-1/2} x^{1-2\sigma}\log(y+2)}{\eta}\right) \right) \\
+ \mathcal{O}\left( \eta t^{2-2\sigma}\right) + \mathcal{O}\left( t^{4/3-2\sigma} \log t\right),
\end{multline*}
where we have used that $\mathcal{E}=\mathcal{O}(1)$. Summing over the error terms we find that for $\sigma<1/2$
\begin{multline*}
\sum_{ x \leq t/2\pi} I_s(s) = |K(s)|^2 \sum_{ x\leq t/2\pi} \sum_{m\leq y-\eta} m^{2\sigma-2} \\
+ \mathcal{O}\left( \frac{t^{1-2\sigma}}{\eta^2}\right) + \mathcal{O}\left( \frac{t^{3/2-2\sigma}\log t}{\eta} \right) + \mathcal{O}\left( \eta t^{2-2\sigma}\right) + \mathcal{O}\left( t^{4/3-2\sigma} \log t\right).
\end{multline*}
Write $\eta=t^{-\delta}$. We should choose $\delta>0$ so that it minimizes the maximum of the exponents $(1+2\delta, 3/2+\delta,2-\delta)$ in the $\eta$-dependent error terms. This is found to be at $\delta = 1/4$ and with this choice of $\eta$ the overall error is dominated by
\[ \mathcal{O}\left( t^{7/4-2\sigma} \log t\right). \]
By performing a similar analysis for $\sigma\leq 1/2$, in which case the $\mathcal{O}(x^{-2\sigma}/\eta^2)$ term gives, after summation, an error at worst of order $\mathcal{O}\left( \log t/\eta^2\right)$, we find that again the overall error is dominated by $\mathcal{O}( t^{7/4 -2\sigma} \log t)$. We conclude that for $\sigma \in (0,1)$
\[ \sum_{ x \leq t/2\pi} I_s(s) = |K(s)|^2 \sum_{ x\leq t/2\pi} \left(\sum_{m\leq y-t^{-1/4}} m^{2\sigma-2}\right) + \mathcal{O}\left( t^{7/4-2\sigma} \log t\right). \]

We write the remaining inner sum as
\[ \sum_{m\leq y-t^{-1/4}} = \sum_{m\leq y} - \sum_{y-t^{-1/4} < m \leq y} \]
and the latter sum is empty if $\|y\|>t^{-1/4}$. Using the previous estimate for $A(t,\eta)$
\[ |K(s)|^2 \sum_{ x\leq t/2\pi}\left(\sum_{y-t^{-1/4} < m \leq y} m^{2\sigma-2}\right) = \mathcal{O}\left( t^{7/4-2\sigma} \right)  \]
which can be absorbed into the previous error term. Hence
\[ \sum_{ x \leq t/2\pi} I_x(s) = |K(s)|^2 \sum_{ x\leq t/2\pi} \sum_{m\leq y} m^{2\sigma-2} + \mathcal{O}\left( t^{7/4-2\sigma} \log t\right). \]
Let us temporarily write $N=t/2\pi$. The double sum is
\[ \sum_{x\leq N} \sum_{m \leq N/x} m^{2\sigma-2}, \]
so applying Dirichlet's hyperbola lemma we find
\begin{align*} &\sum_{x\leq N} \sum_{m \leq N/x} m^{2\sigma-2} \\
&\qquad = \sum_{x\leq \sqrt{N}}\sum_{m\leq N/x} m^{2\sigma-2} + \sum_{m\leq \sqrt{N}}m^{2\sigma-2} \sum_{x\leq N/m} 1 - \sum_{x\leq \sqrt{N}} 1 \sum_{m\leq \sqrt{N}} m^{2\sigma-2} \\
&\qquad = \sum_{x\leq \sqrt{N}}\sum_{m\leq N/x} m^{2\sigma-2} + \sum_{m\leq \sqrt{N}} \floor*{\frac{N}{m}} m^{2\sigma-2} - \floor*{\sqrt{N}} \sum_{m\leq \sqrt{N}} m^{2\sigma-2}
\end{align*}
We consider the cases (i) $\sigma \in (0,1/2)$, (ii) $\sigma \in (1/2,1)$ and (iii) $\sigma=1/2$ separately.

\textbf{Case (i)} Since $\sigma\in (0,1/2)$ the sum $\sum_{m\geq 1} m^{2\sigma-2}$ is convergent. By the Euler-Maclaurin formula
\begin{align*} \sum_{x\leq \sqrt{N}}\sum_{m\leq N/x} m^{2\sigma-2} &=\sum_{x\leq \sqrt{N}} \left( \zeta(2-2\sigma) + \mathcal{O}\left( N^{2\sigma-1} x^{1-2\sigma}\right) \right) \\
&= \floor*{ \sqrt{N}} \zeta(2-2\sigma) + \mathcal{O} \left( N^{2\sigma-1} N^{1-\sigma}\right) \\
&= \sqrt{N} \zeta(2-2\sigma) + \mathcal{O}\left( N^\sigma\right).
\end{align*}
And similarly
\begin{align*}
\sum_{m\leq \sqrt{N}} \floor*{\frac{N}{m}} m^{2\sigma-2} &= \sum_{m\leq \sqrt{N}}\left( \frac{N}{m} + \mathcal{O}(1)\right) m^{2\sigma-2} \\
&= N \left( \zeta(3-2\sigma) + \mathcal{O}\left( N^{\sigma-1}\right) \right) + \mathcal{O}(1) \\
&= N\zeta(3-2\sigma) + \mathcal{O}\left( N^\sigma \right)
\end{align*}
and
\begin{align*}
\floor*{\sqrt{N}} \sum_{m\leq \sqrt{N}} m^{2\sigma-2} &= \floor*{\sqrt{N}} \left( \zeta(2-2\sigma) + \mathcal{O}\left( N^{\sigma-1/2}\right)\right) \\
&= \sqrt{N} \zeta(2-2\sigma) + \mathcal{O}\left(N^{\sigma}\right). 
\end{align*}
So the sum is $N\zeta(3-2\sigma)+ \mathcal{O}\left( N^\sigma\right)$.

\textbf{Case (ii)} Since $\sigma\in (1/2,1)$ the sum $\sum_{m\geq 1} m^{2\sigma-2}$ is divergent. Again using the Euler-Maclaurin formula
\begin{align*} \sum_{x\leq \sqrt{N}}\sum_{m\leq N/x} m^{2\sigma-2} &=\sum_{x\leq \sqrt{N}} \left( \frac{ N^{2\sigma-1}x^{1-2\sigma}}{2\sigma-1} + \mathcal{O}\left( N^{2\sigma-2} x^{2-2\sigma} \right)\right) \\
&= \frac{ N^{2\sigma-1}\floor*{\sqrt{N}}^{2-2\sigma}}{(2\sigma-1)(2\sigma-2)} + \mathcal{O} \left( N^{2\sigma-2} N^{3/2-\sigma}\right) \\
&= \frac{ N^{\sigma}}{(2\sigma-1)(2\sigma-2)} + \mathcal{O} \left( N^{\sigma-1/2}\right).
\end{align*}
And similarly
\begin{align*}
\sum_{m\leq \sqrt{N}} \floor*{\frac{N}{m}} m^{2\sigma-2} &= \sum_{m\leq \sqrt{N}}\left( \frac{N}{m} + \mathcal{O}(1)\right) m^{2\sigma-2} \\
&= N \left( \zeta(3-2\sigma) + \mathcal{O}\left( N^{\sigma-1}\right) \right) + \mathcal{O}(N^{\sigma-1/2}) \\
&= N\zeta(3-2\sigma) + \mathcal{O}\left( N^\sigma \right)
\end{align*}
and
\begin{align*}
\floor*{\sqrt{N}} \sum_{m\leq \sqrt{N}} m^{2\sigma-2} &= \floor*{\sqrt{N}} \left( \frac{\floor*{\sqrt{N}}^{2\sigma-1}}{2\sigma-1} + \mathcal{O} \left( N^{\sigma-1}\right)\right) \\
&= \frac{N^\sigma}{2\sigma-1} + \mathcal{O}\left(N^{\sigma-1/2}\right). 
\end{align*}
So again the sum is $N\zeta(3-2\sigma)+ \mathcal{O}\left( N^\sigma\right)$.

\textbf{Case (iii)} For $\sigma=1/2$ we use Euler Maclaurin again, along with Stirling's estimate
\begin{align*} \sum_{x\leq \sqrt{N}}\sum_{m\leq N/x} m^{-1} &=\sum_{x\leq \sqrt{N}} \left( \log N - \log x + \mathcal{O}(1) \right) \\
&= \floor*{\sqrt{N}} \log N - \left( \floor*{\sqrt{N}} \log \floor*{\sqrt{N}} - \log \floor*{\sqrt{N}} \right) + \mathcal{O}\left( \sqrt{N}\right) \\
&= \tfrac{1}{2} \sqrt{N} \log N + \mathcal{O}\left( \sqrt{N}\right).
\end{align*}
And similarly
\begin{align*}
\sum_{m\leq \sqrt{N}} \floor*{\frac{N}{m}} m^{-1} &= \sum_{m\leq \sqrt{N}}\left( \frac{N}{m} + \mathcal{O}(1)\right) m^{-1} \\
&= N \left( \zeta(2) + \mathcal{O}\left( N^{-1/2}\right) \right) + \mathcal{O}\left( \log N\right) \\
&= N\zeta(2) + \mathcal{O}\left( \sqrt{N} \right)
\end{align*}
and
\begin{align*}
\floor*{\sqrt{N}} \sum_{m\leq \sqrt{N}} m^{-1} &= \floor*{\sqrt{N}} \left( \log \sqrt{N} + \mathcal{O}(1) \right) \\
&= \tfrac{1}{2}\sqrt{N} \log N + \mathcal{O}\left( \sqrt{N}\right). 
\end{align*}
So the sum is $N\zeta(2)+ \mathcal{O}\left(\sqrt{N}\right)$.

Collecting these results together, we conclude that for $\sigma\in (0,1)$
\[ |K(s)|^2 \sum_{x\leq t/2\pi} \sum_{m\leq y} m^{2\sigma-2} = \left(\frac{t}{2\pi}\right)^{2-2\sigma} \zeta(3-2\sigma) + \mathcal{O}\left( t^{1-\sigma} \right), \]
and using the previous estimates we arrive at Theorem \ref{thm3}.

\end{document}